\pdfoutput=1
\documentclass{article}

\usepackage[a4paper]{geometry}

\usepackage{amsmath}
\usepackage{amssymb}
\usepackage{amsthm}
\usepackage[utf8]{inputenc} 
\usepackage[T1]{fontenc} 
\usepackage{xcolor}
\usepackage{subcaption}
\usepackage{booktabs}
\usepackage{natbib}
\usepackage{lmodern}
\usepackage{orcidlink}

\usepackage{algorithm}
\usepackage{algorithmic}

\newtheorem{remark}{Remark}

\usepackage{hyperref}

\newcommand{\revised}[1]{#1}

\newcommand\BibTeX{{\rmfamily B\kern-.05em \textsc{i\kern-.025em b}\kern-.08em
T\kern-.1667em\lower.7ex\hbox{E}\kern-.125emX}}

\begin{document}

\title{Performance of linear solvers in tensor-train format on current multicore architectures}
\author{
	Melven Röhrig-Zöllner\orcidlink{0000-0001-9851-5886}\thanks{German Aerospace Center (DLR), Institute of Software Technology} \and
	Manuel Becklas\footnotemark[1] \and
	Jonas Thies\orcidlink{0000-0001-9231-9999}\thanks{Delft University of Technology, Delft Institute of Applied Mathematics} \and
	Achim Basermann\orcidlink{0000-0003-3637-3231}\footnotemark[1] }



\maketitle

\begin{abstract}
\revised{
	Tensor networks are a class of algorithms aimed at reducing the computational complexity of high-dimensional problems.
	They are used in an increasing number of applications, from quantum simulations to machine learning.
	Exploiting data parallelism in these algorithms is key to using modern hardware.
	However, there are several ways to map required tensor operations onto linear algebra routines (``building blocks'').
	Optimizing this mapping impacts the numerical behavior, so computational and numerical aspects must be considered hand-in-hand.
	In this paper we discuss the performance of solvers for low-rank linear systems in the tensor-train format (also known as matrix-product states).
	We consider three popular algorithms: TT-GMRES, MALS, and AMEn.
	We illustrate their computational complexity based on the example of discretizing a simple high-dimensional PDE in, e.g., $50^{10}$ grid points.
	This shows that the projection to smaller sub-problems for MALS and AMEn reduces the number of floating-point operations by orders of magnitude.
	We suggest optimizations regarding orthogonalization steps, singular value decompositions, and tensor contractions.
	In addition, we propose a generic preconditioner based on a TT-rank-1 approximation of the linear operator.
	Overall, we obtain roughly a 5x speedup over the reference algorithm for the fastest method (AMEn) on a current multicore CPU.
}
\end{abstract}




\clearpage

\section{Introduction}
Low-rank tensor methods provide a way to approximately solve problems that would otherwise require huge amounts of memory and computing time.
Many ideas in this field arise from quantum physics. E.g., the global state of a quantum system with $N$ two-state particles can be expressed as a tensor of dimension $2^N$.
In this setting, interesting states are for example given by the eigenvectors of the smallest eigenvalues of the Hamiltonian of the system--a Hermitian linear operator that describes the energy of the system.
Therefore, most work focuses on solving eigenvalue problems for Hermitian/symmetric operators using the DMRG method \citep{Wilson1983renormalization_group,White1992,Schollwoeck2005}.
However, linear systems in low-rank tensor formats also arise in interesting applications for example for solving high-dimensional or parameterized partial differential equations, see, e.g., \cite{Kressner2011,Dahmen2015,Dolgov2019}.
In addition, linear solvers and eigenvalue solvers are closely related and many successful methods for finding eigenvalues are based on successive linear solves.
This paper addresses iterative methods for solving linear systems (symmetric and non-symmetric) in the tensor-train (TT) format for the case where the individual dimensions are not tiny, i.e., for systems of dimension $n^d\times n^d$ with $n\gg 2$.
We employ the TT format (called matrix-product states in physics) as it is a simple and common low-rank tensor format.
Most of the ideas are transferable to other low-rank tensor formats (at least to loop-free tensor networks).
Our work considers the TT-GMRES algorithm \citep{Dolgov2013a,Ballani2012}, the modified alternating least-squares (MALS) algorithm (DMRG for linear systems) \citep{Holtz2012,Oseledets2012}, and the Alternating Minimal Energy (AMEn) algorithm \citep{Dolgov2014a}.
We show numerical improvements and performance improvements of the underlying operations and focus on a single CPU node.
These improvements are orthogonal to the parallelization for distributed memory systems presented in \cite{Daas2022}, so the suggestions from both \cite{Daas2022} and this paper could be combined in the future.
As the resulting complete linear solver requires a tight interplay of different algorithmic components, we discuss the behavior of the different numerical methods involved for the TT format.
An alternative class of methods to solve linear systems in TT format consists of Riemannian optimization on the manifold of fixed-rank tensor-trains, see \cite{Kressner2016}.
This results in a nonlinear optimization problem and is therefore not in the scope of this paper. However, it partly requires similar underlying operations.

The paper is organized as follows:
First, we start with required numerical background in \autoref{sec:background_and_notation} and also introduce relevant performance metrics for today's multicore computers.
Then in \autoref{sec:numerical_algorithms}, we discuss the involved high-level numerical algorithms TT-GMRES, MALS, and AMEn.
\revised{Based on an example, we illustrate the numerical behavior of these algorithms and compare their complexity in \autoref{sec:comparison_of_methods}.
Afterwards, in \autoref{sec:building_blocks}, we analyze and optimize the underlying building blocks.}
We conclude in \autoref{sec:conclusion} with a short summary and open points for future work.

\section{Background and notation}\label{sec:background_and_notation} 
In this section, we provide the required background concerning numerics and performance.

\subsection{Numerical background}
\revised{%
We first introduce required matrix decompositions and a notation for the considered algorithms.
}%

\subsubsection{Matrix decompositions}
As matrix decompositions are heavily used as steps in tensor-train algorithms, we repeat some basic properties of QR and SVD decompositions from the literature, see e.g., \cite{Golub2013, Higham2002}.

A QR-decomposition is a factorization of a matrix $M\in\mathbf R^{n_1\times n_2}$ with column rank $r$ into a matrix $Q$ with orthonormal columns and an upper triangular part $R$:
\begin{align}
	M &= QR, &&\text{with}&& Q^TQ=I,\; Q\in\mathbf R^{n_1 \times r},\; R\in\mathbf R^{r\times n_2}.
\end{align}
For (numerically) rank-deficient $M$, one can employ a pivoted QR-decomposition
\begin{align}
MP &= QR,&&\Leftrightarrow& M &= QRP^{-1},
\end{align}
where $P$ is a permutation matrix.
The pivoted QR-decomposition can be computed in a numerically robust way.
However, it cannot (safely) be used to approximate $M$ with a lower rank matrix based on the size of the pivot elements (diagonal entries of $R$) by $Q_{:,1:r'}R_{1:r',:}$, $r'<r$ as the worst case error grows with $O(2^r)$, see \cite{Higham1990,Kawamura2021}.

The singular value decomposition (SVD) in contrast provides the best approximation of lower rank:
\begin{align}
M &= USV^T&&\Rightarrow& \|M - U_{:,1:r'} S_{1:r',1:r'} V_{:,1:r'}^T\|_F = \min_{\tilde M, \operatorname{rank}(\tilde M)=r'}\| M - \tilde M \|_F,
\end{align}
where $U\in\mathbf R^{n_1\times r}$, $V\in\mathbf R^{n_2\times r}$ are the matrices of the orthonormal left/right singular vectors and $S=\operatorname{diag}(s_1, s_2, \dots, s_r)$ is composed of the singular values $s_1 \ge s_2 \ge \dots \ge s_r > 0$.

\subsubsection{Tensor-train decomposition}
In higher dimensions, there is no unique way to decompose a tensor into factors and to define its rank(s) (variants are e.g., the Tucker and the CANDECOMP/PARAFAC (CP) decompositions, see the review \cite{Kolda2009}).
We focus on the tensor-train format which decomposes a tensor $X\in\mathbf R^{n_1 \times n_2 \times \cdots \times n_d}$ into $d$ three-dimensional sub-tensors $X_1, \dots, X_d$:
\begin{align}
X &= X_1 \Join X_2 \Join \cdots \Join X_d && \text{with}& X_k \in\mathbf R^{r_{k-1}\times n_k \times r_k}, r_0 = r_d = 1.
\end{align}
Here, $(\cdot \Join \cdot)$ is the contraction of the last dim.\ of the left operand with the first dim.\ of the right operand:
\begin{align}
X_k \Join X_{k+1} &= \sum_i (X_k)_{:,:,i} (X_{k+1})_{i,:,:} \; \in \mathbf R^{r_{k-1} \times n_k \times n_{k+1} \times r_{k+1}}.
\end{align}
The TT decomposition of a given tensor is not unique: it is invariant with respect to multiplying one sub-tensor by a matrix and the next with its inverse.
\revised{More precisely, for $M,N^T\in\mathbf R^{r_k \times r_k'}$ with $MN=I$:}
\begin{align}
\bar X_k \Join \bar X_{k+1} &= X_k \Join X_{k+1} &&\text{for}& \bar X_k := X_k\Join M, \quad \bar X_{k+1} := N \Join X_{k+1}.
\end{align}
The smallest possible dimensions $(r_1, \dots r_{d-1})$ that allow to represent $X$ denote the TT ranks of $X$ with the maximal rank $r:=\operatorname{rank}(X)=\max(r_1, \dots, r_{d-1})$.
If $X$ has rank-1 in the TT format, we can also write it as a generalized dyadic product of a set of vectors:
\begin{align*}
	X &= (X_1)_{1,:,1} \otimes (X_2)_{1,:,1} \otimes \cdots \otimes (X_d)_{1,:,1}.
\end{align*}

\subsubsection{Tensor unfolding and orthogonalities}
We define a general \emph{reshape} operation to reinterpret the entries of a tensor as a tensor of different dimensions:
\begin{align*}
	\operatorname{reshape}\left(X, \begin{pmatrix} \bar n_1 & \cdots & \bar n_{\bar d} \end{pmatrix}\right) &:= \bar X \in\mathbf R^{\bar n_1\times\cdots\times \bar n_{\bar d}} \quad\text{with}\\
	(\bar X)_{\bar i_1,\dots,\bar i_{\bar d}} = (X)_{i_1, \dots, i_d} \quad&\text{for}\quad \bar i_1 + \bar i_2 \bar n_1 + \bar i_3 \bar n_1 \bar n_2 + \dots = i_1 + i_2 n_1 + i_3 n_1 n_2 + \cdots .
\end{align*}
With this, we define the \emph{left-unfolding} that combines two dimensions of a 3d tensor $X_k$ to obtain a matrix:
\begin{align}
(X_k)_\text{left} &:= \operatorname{reshape}\left(X_k, \begin{pmatrix}r_{k-1}n_k& r_k\end{pmatrix}\right) \in \mathbf R^{r_{k-1}n_k \times r_k}.
\end{align}
Similarly, we define the \emph{right-unfolding}:
\begin{align}
(X_k)_\text{right} &:= \operatorname{reshape}\left(X_k, \begin{pmatrix}r_{k-1} & n_k r_k\end{pmatrix}\right) \in \mathbf R^{r_{k-1} \times  n_k r_k}.
\end{align}
We denote a 3d tensor $X_k$ as \emph{left-orthogonal} if the columns of the left-unfolding are orthonormal ($(X_k)_\text{left}^T (X_k)_\text{left} = I$) and \emph{right-orthogonal} if the rows of the right-unfolding are orthonormal ($(X_k)_\text{right} (X_k)_\text{right}^T=I$).
From the TT format, we can build an SVD of an unfolding of $X\in\mathbf R^{n_1\times n_2 \cdots \times n_d}$ into a set of left and right dimensions
\begin{align}
USV^T &= (X)_{\text{unfold}_j} := \operatorname{reshape}\left(X, \begin{pmatrix}n_1\cdots n_j & n_{j+1} \cdots n_d\end{pmatrix}\right)
\end{align}
by left-orthogonalizing ($X_1, \dots, X_j$) and right-orthogonalizing ($X_{j+1}, \dots, X_d$):
\begin{align}
X &= \underbrace{\bar X_1 \Join \cdots \Join \bar X_j}_{U} \Join S \Join \underbrace{\bar X_{j+1} \Join \cdots \Join \bar X_d}_{V^T}.
\end{align}

\subsubsection{Tensor-train vectors and operators}
A tensor-train operator is a tensor in TT format where we combine pairs of dimensions $(n_i\times m_i)$ with the form $\mathcal A\in\mathbf R^{(n_1\times m_1) \times (n_2 \times m_2) \times \cdots \times (n_d \times m_d)}$
\begin{align}
\mathcal A &= A_1 \Join A_2 \Join \cdots A_d&&\text{with}& A_k\in\mathbf R^{r^A_{k-1}\times n_k\times m_k\times r^A_k}, r^A_0=r^A_d=1
\end{align}
which defines a linear mapping: $\mathcal A: \mathbf R^{m_1 \times m_2 \times \cdots \times m_d} \to \mathbf R^{n_1 \times n_2 \times \cdots \times n_d}$.
For the scope of this paper, we only consider quadratic regular operators ($m_k = n_k$).
Again, the sub-tensor dimensions $(r^A_1, \dots r^A_{d-1})$ denote the TT ranks of the operator with a maximal rank of $r^A:=\operatorname{rank}(\mathcal A)=\max(r^A_1, \dots r^A_{d-1})$.
This definition allows an efficient application of a TT operator on a TT ``vector'' directly in the TT format:
\begin{align}
\mathcal A X &= Y = Y_1 \Join Y_2 \Join \cdots  \Join Y_n\qquad \text{with}\nonumber\\
Y_k &:= \operatorname{reshape}\left( \sum_i (A_k)_{:,:,i,:} (X_k)_{:,i,:}, \begin{pmatrix} r^A_{k-1} r_{k-1} & n_k & r^A_k r_k \end{pmatrix} \right).
\end{align}
The resulting TT decomposition has ranks $(r^A_1 r_1, \dots, r^A_{d-1} r_{d-1})$ and thus $\operatorname{rank}(Y) \le r^A r$.

With all these definitions at hand, we can specify the main problem considered in this paper:
Given a low-rank linear operator in TT format $\mathcal A_\text{TT}:=(A_1, \dots A_d)$, a low-rank right-hand side (RHS) $B_\text{TT}:=(B_1, \dots, B_d)$,
and a desired residual tolerance $\epsilon_\text{tol}$, find an approximate low-rank solution $X_\text{TT}:=(X_1, \dots, X_d)$ with
\begin{align}\label{eq:linear_system}
\left\|\mathcal A_\text{TT} X_\text{TT} - B_\text{TT}\right\|_F &\le \epsilon_\text{tol}.
\end{align}

\subsection{Performance characteristics on today's multicore CPU systems}
Today's hardware features multiple levels of parallelism and memory that we need to exploit to efficiently use the available compute capacity:
A supercomputer is composed of a number of nodes connected by a network (distributed memory parallelism).
Each node contains one or more multicore CPUs with access to a main memory (shared memory parallelism).
However, the different CPUs access different memory domains with different speed (NUMA architecture): Usually the access to the memory banks directly connected to one CPU is faster.
Each CPU ``socket'' consists of multiple cores ($\lesssim 100$ in 2023) with a hierarchy of caches. 
Inside one core, SIMD units perform identical calculations on a small vector of floating-point numbers. 
We focus on the performance on a single multicore CPU, but the algorithms considered can also run in parallel on a cluster (see \cite{Daas2022}).
Many supercomputers nowadays use GPUs which is not discussed further in this paper.

\subsubsection{Roofline performance model}
To obtain a simpler abstraction of the hardware, we employ the Roofline \citep{Williams2009} performance model.
\revised{The Roofline model distinguishes between computations (floating-point operations) and data transfers:}
The maximal performance one can achieve on a given hardware is $P_\text{peak}$ [GFlop/s].
The bandwidth of data transfers from main memory is $b_s$ [GByte/s]. If all data fits into some cache level, the appropriate cache bandwidth is used instead.
These are the required hardware characteristics.
The considered characteristics of one building block of an algorithm (e.g., one nested loop) are the number of required floating-point operations $n_\text{flops}$ and the volume of the data transfers $V_\text{read/write}$.
Their ratio is called computational intensity $I_c := n_\text{flops} / V_\text{read/write}$ [Flop/Byte].
Assuming that the data transfers and operations overlap perfectly, this results in the following performance:
\begin{align}
P_\text{roofline} &= \min\left(P_\text{peak}, I_c b_s\right).
\end{align}
If the compute intensity is low ($I_c \ll P_\text{peak}/b_s$), the building block is \emph{memory bound}.
If in contrast the compute intensity is high ($I_c \gg P_\text{peak}/b_s$), the building block is \emph{compute bound}.
We specifically split the algorithms considered in this paper into smaller parts (``building blocks'') because they feature not one dominating operation but are composed of multiple different blocks with different performance characteristics.

\subsubsection{Memory and cache performance}
In addition to the model above, some details of the memory hierarchy play a crucial role for the algorithms at hand (see \cite{Hager2010} for more details):
First, modifying memory is often slower than reading it. In order to write to main memory, the memory region is usually first transferred to the cache, modified there and written back (write-allocate).
A special CPU instruction allows to avoid this and directly stream to memory (non-temporal store).
A common technique to improve the performance of data transfers is to avoid writing large (temporary) data when the algorithm can be reformulated accordingly (write-avoiding), see \cite{Carson2016}.
Second, the CPU caches are organized in cache lines: This means that, e.g., 8 double precision values are transferred together, always starting from a memory address divisible by the cache line size. So the data locality--e.g., which index is stored contiguously--is important.

In addition, today's CPUs use set-associative caches that allow the mapping of one memory address to a fixed set of cache lines.
Due to this, memory addresses with a specific distance (e.g., $1024$ double numbers) are mapped to the same cache set and the cache effectiveness is dramatically reduced when data is accessed with specific ``bad'' strides (\emph{cache thrashing}).
This easily occurs for tensor operations if the product of some dimensions is close to a power of two.
A common solution for operations on 2d arrays is padding: adding a few ignored zero rows in a matrix such that the stride is at least a few cache lines bigger than some power of two.
This becomes more complicated in higher dimensions as one either obtains a complex indexing scheme or one needs to perform calculations with zeros.
This is discussed in more detail in \autoref{sec:faster_contractions}.

\section{Numerical algorithms}\label{sec:numerical_algorithms}
In this section, we discuss three different methods to approximately solve a linear system in TT format as in \autoref{eq:linear_system}.
We start with a general purpose (``global'') approach based on Krylov subspace methods, TT-GMRES \citep{Dolgov2013a}, and present some improvements for the TT format.
Then, we consider the more specialized (``local'') MALS \citep{Holtz2012,Oseledets2012}, which optimizes pairs of sub-tensors similar to DMRG.
Afterwards, we discuss the (``more local'') AMEn \citep{Dolgov2014a} method, which iterates on one sub-tensor after another.
Finally, we present a simple yet effective preconditioner in TT format to accelerate convergence of TT-GMRES.

\subsection{Krylov methods: TT-GMRES}
All methods that apply the linear operator $\mathcal A_\text{TT}$ on linear combinations of previously calculated directions produce solutions from the Krylov subspace $\mathbf K_k(\mathcal A_\text{TT}; V_\text{TT}):=\operatorname{span}\left\{V_\text{TT}, \mathcal A_\text{TT} V_\text{TT}, \dots, \mathcal A_\text{TT}^{k-1} V_\text{TT}\right\}$ where $V_\text{TT}$ is usually the initial residual of the problem.
Different Krylov subspace methods then select the ``best'' solution from the subspace $\mathbf K_k$ according to different definitions of ``best'', see \cite{Vorst2003} and \cite{Saad2003} for a detailed discussion.
As we usually only approximate intermediate steps in TT arithmetic, we effectively employ inexact Krylov methods which are discussed thoroughly in \cite{Simoncini2003,Eshof2004}.
In this paper, we consider the TT-GMRES method \citep{Dolgov2013a} for non-symmetric problems. For symmetric operators $\mathcal A_\text{TT}$, we simply omit unneeded steps to obtain a MINRES variant.
However, all considerations here effectively hold for other Krylov subspace methods as well.

\subsubsection{Arithmetic operations in tensor-train format}
Krylov methods require the following operations which can be performed directly in the TT format (all introduced in \cite{Oseledets2011}):
Applying the operator to a vector ($Y_\text{TT}\gets \mathcal A_\text{TT} X_\text{TT}$), dot products ($\alpha \gets \langle X_\text{TT}, Y_\text{TT} \rangle$) and scaled additions ($Y_\text{TT} \gets \alpha X_\text{TT} + Y_\text{TT}$) of two tensor-trains.
To reduce the computational complexity, TT truncation ($\tilde X_\text{TT} \gets \operatorname{trunc}_\delta(X_\text{TT})$) approximates a tensor-train with another tensor-train with lower rank (see TT-rounding algorithm in \cite{Oseledets2011}):
\begin{align*}
	\left\|X_\text{TT} - \operatorname{trunc}_{\delta}(X_\text{TT})\right\|_F &\le \delta .
\end{align*}
With these operations, we can perform a variant of the GMRES algorithm with additional truncation operations, see \autoref{alg:TT-GMRES}.
\revised{%
This idea was first discussed in \cite{Ballani2012} for the more general H-Tucker format with a slightly different projection.
We employ a variant based on \cite{Dolgov2013a}. The numerical stability of TT-GMRES is analyzed in more detail in \cite{Coulaud2022}.
}%
\begin{algorithm}
\caption{TT-GMRES with modified Gram-Schmidt}
\label{alg:TT-GMRES}
\begin{algorithmic}[1]
\REQUIRE{Linear operator $\mathcal A_\text{TT}: \mathbf R^{n_1 \times \cdots \times n_d} \to \mathbf R^{n_1 \times \cdots \times n_d}$ and RHS $B_\text{TT} \in \mathbf R^{n_1 \times \cdots \times n_d}$,\\desired tolerance $\epsilon$, max. number of iterations $m$, estimated condition number $c$}
\ENSURE{Approximate solution $X_\text{TT}$ with $\left\|B_\text{TT} - \mathcal A_\text{TT} X_\text{TT}\right\|_F \lesssim \epsilon \left\|B_\text{TT}\right\|_F$}
	\STATE {$\gamma_0 \gets \| B_\text{TT} \|_F$}
	\STATE {$V_{TT,1} \gets 1/\gamma_0 B_\text{TT}$}
	\FOR {$i=1, \dots, m$}
		\STATE {Choose tolerance $\delta_i = \frac{0.5\epsilon}{cm}\frac{\gamma_0}{\gamma_{i-1}}$} \hfill\COMMENT{In \cite{Dolgov2013a}: $\delta_i = \epsilon\frac{\gamma_0}{\gamma_{i-1}}$}
		\STATE {$W_\text{TT} \gets \operatorname{trunc}_{0.5\delta_i/(i+1)}(\mathcal A_\text{TT} V_{\text{TT},i})$}
		\FOR {$j=1, \dots, i$}
			\STATE {$h_{i,j} \gets \langle V_{\text{TT},j}, W_\text{TT} \rangle $} \hfill\COMMENT{MINRES: $h_{i,j}=0, j < i-1$}
			\STATE {$W_\text{TT} \gets \operatorname{trunc}_{0.5\delta_i/(i+1)}( W_\text{TT} - h_{i,j} V_{\text{TT},j} )$}
		\ENDFOR
		\STATE {$W_\text{TT} \gets \operatorname{trunc}_{0.5\delta_i}(W_\text{TT})$}
		\STATE {$h_{i+1,i} \gets \| W_\text{TT} \|_F$}
		\STATE {$V_{\text{TT},i+1} \gets 1/h_{i+1,i} W_\text{TT}$}
		\STATE {$y \gets \arg\min_y \|H y - \gamma_0 e_1\|_2$ for $H = (h_{k,l}),\;k=1,\dots,i+1,\;l=1,\dots,i$}
		\STATE {$\gamma_i \gets \|H y - \gamma_0 e_1\|_2$}
		\IF {$\gamma_i / \gamma_0 \le 0.5\epsilon$}
			\STATE {$X_\text{TT} \gets 0$}
			\FOR {$j=1, \dots, i$}
				\STATE {$X_\text{TT} \gets \operatorname{trunc}_{0.5\epsilon/(ci)} (X_\text{TT} + y_j V_{\text{TT},j})$}
			\ENDFOR
			\RETURN
		\ENDIF
	\ENDFOR
	\STATE {Abort: not converged in $m$ steps!}
\end{algorithmic}
\end{algorithm}
We remark that we use a more strict truncation tolerance than suggested in \cite{Dolgov2013a} based on the analysis of inexact Krylov methods in \cite{Simoncini2003}.
However, in \cite{Simoncini2003} only inaccurate applications of the linear operator are considered (as in line~5 of \autoref{alg:TT-GMRES}).
We also truncate in each step of the orthogonalization (line~8) and once afterwards (line~10).
We can still express the error as an error in the operator of the form (see eq.~(2.2) in \cite{Simoncini2003}):
\begin{align*}
	\left(A+E_i\right)v_i &= V_{i+1} H_i, \qquad i=1, \dots, m .
\end{align*}
We denote the errors of all truncation operations in one Arnoldi iteration with $\Delta w^{(0)}$ (line~5), $\Delta w^{(j)}$ (line~8) and $\Delta w^{(i+1)}$.
Then, we obtain for the error:
\begin{align*}
	E_iv_i &= \Delta w^{(i+1)} + \sum_{j=0}^i\left(I - \sum_{k=j+1}^i v_k v_k^T\right)\Delta w^{(j)} .
\end{align*}
For orthogonal basis vectors $v_k$ of the Krylov subspace, the error of the Arnoldi iteration is bounded by:
\begin{align*}
	\|E_i\| 
	        &\le \|\Delta w^{(i+1)}\|_F + \sum_{j=0}^i\|\Delta w^{(j)}\|_F
			\le 0.5\delta_i\|w^{(i+1)}\|_F + \sum_{j=0}^i \frac{0.5\delta_i}{i+1}\|w^{(j)}\|_F \le \delta_i \|A\| .
\end{align*}
Of course, due to the truncations, one easily looses the orthogonality of the basis vectors $v_k$, see discussion in \autoref{subsubsec:SIMGS} below.
Suitable tolerances $\delta_i$ require the condition number $c=\kappa(H_m)$ that we estimate using the parameter $c\approx\kappa(A)$ as suggested in \cite{Simoncini2003}.
So we obtain the following bound for the difference between the true residual vector $r_*$ and the inexact residual vector $\tilde r_*$ (see equation~(5.8) in \cite{Simoncini2003}):
\begin{align*}
	\|r_* - \tilde r_*\|_F &\le 0.5\epsilon .
\end{align*}
The factor $0.5$ ensures that the true residual norm is smaller than the desired tolerance:
\begin{align*}
	\|r_*\| &= \|r_* - \tilde r_* + \tilde r_*\| \le \| r_* - \tilde r_*\| + \|\tilde r_*\| \le \epsilon .
\end{align*}
In our experiments, we use an optimized variant (see \autoref{sec:building_blocks_svd_and_qr}) of the standard TT truncation algorithm.
An alternative randomized truncation algorithm is presented in \cite{Daas2023} for truncating sums of multiple tensor-trains (e.g., only truncating in line~10 of \autoref{alg:TT-GMRES} and not in line~8).

\subsubsection{Improved Gram-Schmidt orthogonalization}\label{subsubsec:SIMGS}
Above, we assumed that the resulting Krylov basis vectors are orthogonal.
However, as the modified Gram-Schmidt orthogonalization is only applied approximately (truncations in line~8 and~10), this assumption is usually violated.
As a result, the true residual norm might not be smaller than the prescribed tolerance $\epsilon$ even though the approximate residuals converge.
To compensate for the inaccurate orthogonalization, one can prescribe a more strict truncation tolerance as discussed in section~6 of \cite{Simoncini2003}.
Another common approach consists in re-orthogonalization:
We employ the following specialized variant of a modified Gram-Schmidt iteration.

In the TT format calculating a scalar product is much faster than a truncated scaled addition (axpby) as discussed in more detail in \autoref{sec:building_blocks}.
So we can perform additional scalar products to reorder Gram-Schmidt iterations and perform selective re-orthogonalization as shown in \autoref{alg:SIMGS} to increase the robustness.
\begin{algorithm}
\caption{Selective iterated modified Gram-Schmidt (SIMGS) in TT format}
\label{alg:SIMGS}
\begin{algorithmic}[1]
	\REQUIRE{Orthonormal previous directions $V_{\text{TT},j}$, $j = 1, \dots, i$, new direction $W_\text{TT}$, \\
	         tolerance $\delta_i$, max.\ re-orthogonlization iterations $k_\text{max}$}
	\ENSURE{New normalized direction $V_{\text{TT},i+1}$ with $W_\text{TT} \approx \sum_{j=1}^{i+1} h_j V_{\text{TT},j}$,\\
			and $\left|\langle V_{\text{TT},i+1} , V_{\text{TT},j} \rangle\right| \lesssim \delta_i$, $j=1, \dots, i$}
	\STATE {$W_\text{TT} \gets \operatorname{trunc}_{0.5\delta_i/(2i+1)}\left(W_\text{TT}\right)$}
    \STATE {$h_j \gets 0$, $j = 1, \dots i$}
	\FOR {$k=1, \dots, k_\text{max}$}
		\STATE {Calculate $g_j \gets \langle V_{\text{TT},j}, W_\text{TT} \rangle / \|W_\text{TT}\|_F$, $j = 1, \dots i$}
		\STATE {\textbf{if} $\|g\|_\infty \le \delta_i$ \textbf{break}}
		\FOR {$j = \arg\max_l |g_l|$ and $|g_j| > \delta_i$}
			\STATE {$\beta \gets \langle V_{\text{TT},j}, W_\text{TT} \rangle$}
			\STATE {$W_\text{TT} \gets \operatorname{trunc}_{0.5\delta_i/(2i+1)}( W_\text{TT} - \beta V_{\text{TT},j} )$}
			\STATE {$h_j \gets h_j + \beta, \qquad g_j \gets 0$}
		\ENDFOR
	\ENDFOR
	\STATE {$W_\text{TT} \gets \operatorname{trunc}_{0.5\delta_i}\left(W_\text{TT}\right)$}
	\STATE {$h_{i+1} \gets \|W_\text{TT}\|_F, \qquad V_{\text{TT},i+1} \gets W_\text{TT} / h_{i+1}$}
\end{algorithmic}
\end{algorithm}
This omits subtracting directions that are already almost orthogonal in order to avoid growing the TT-ranks.
See \cite{Leon2012} and the references therein for a detailed discussion on different Gram-Schmidt orthogonalization schemes.
Here, we again need to use sufficiently small truncation tolerances to fulfill the requirements of the outer inexact GMRES method.
The factor $2i$ is an estimate for the number of inner iterations (line~9) as usually orthogonalizing ``twice is enough'' \citep{Parlett1998,Giraud2005}.
And we choose $k_\text{max}=4$ in all our experiments as this was sufficient for the cases we investigated.

\subsubsection{Tensor-train ranks for problems with a displacement structure}
Even with truncations after each tensor-train addition, the tensor-train ranks can grow exponentially in the worst case:
\begin{align*}
	\operatorname{rank}(V_{\text{TT},i+1}) &\le \operatorname{rank}(\mathcal A_\text{TT}) \operatorname{rank}(V_{\text{TT},i}) + \sum_{j=0}^{i} \operatorname{rank}(V_{\text{TT},j}) .
\end{align*}
We observe only a much smaller growth for some special linear operators $\mathcal A_\text{TT}$.
In particular, we consider linear operators with a displacement/Laplace structure:
\begin{align*}
	(\mathcal A X)_{\text{unfold}_j} &= A_{j,\text{left}} X_{\text{unfold}_j} + X_{\text{unfold}_j} A_{j,\text{right}}, \qquad j = 1, \dots, d-1 .
\end{align*}
For those operators, the rank of the solution is bounded if the right-hand side also has small rank as discussed for several tensor formats in \cite{Shi2021}.
We obtain the following expression for applying the operator $k$ times:
\begin{align*}
	(\mathcal A^k X)_{\text{unfold}_j} &= \sum_{i=0}^k \binom{k}{i} A_{j,\text{left}}^{k-i} X_{\text{unfold}_j} A_{j,\text{right}}^{i} .
\end{align*}
%
Similarly, for any matrix polynomial $p_k$ of degree $k$, we get
\begin{align*}
	(p_k(\mathcal A) X)_{\text{unfold}_j} &= \sum_{i=0}^k \bar p_{k-i}(A_{j,\text{left}}) X_{\text{unfold}_j} \hat p_i(A_{j,\text{right}})
\end{align*}
with appropriate sets of polynomials $\bar p_i$ and $\hat p_i$.
As the $j$th TT rank is just the rank of the $j$th unfolding, this results in at most a linear growth of ranks of Krylov subspace basis vectors:
\begin{align}
	\operatorname{rank}(V_{\text{TT},i}) &\le (i+1) \operatorname{rank}(V_{\text{TT},0}).
\end{align}
However, we will see in \autoref{subsec:behavior_of_tt_ranks} that this only holds in exact arithmetic.
If we do not calculate the Krylov basis accurately enough, the TT ranks might again grow exponentially.

\subsection{Modified Alternating Least-Squares (MALS)}
Krylov methods like TT-GMRES consider the linear operator as a black box.
However, we can also exploit the tensor-train structure of the problem and project it onto the subspace of one or several sub-tensors.
This is the idea of the ALS and MALS methods discussed in \cite{Holtz2012}.
In principle, MALS is identical to the famous DMRG method \citep{White1992,Schollwoeck2005} for eigenvalue problems from quantum physics.
Here, we describe it from the point of view of numerical linear algebra.
Using all but two sub-tensors of the current approximate solution, we define the operator:
\begin{align}\label{eq:mals_projection_operator} 
	\mathcal V_{j,2}&: \mathbf R^{r_{j-1}\times n_j \times n_{j+1} \times r_{j+1}} \to \mathbf R^{n_1 \times \cdots \times n_d}, \quad \text{with} \nonumber\\
	\mathcal V_{j,2}Y&= X_1\Join \cdots \Join X_{j-1} \Join Y \Join X_{j+2} \Join \cdots \Join X_d.
\end{align}
If sub-tensors ($X_1, \dots X_{j-1}$) are left-orthogonal and sub-tensors ($X_{j+2},\dots X_d$) are right-orthogonal, the operator $\mathcal V_{j,2}$ has orthonormal columns: $\mathcal V_{j,2}^T \mathcal V_{j,2} = I$.
Using $\mathcal V_{j,2}$, we can project the problem onto the subspace of the sub-tensors $\{X_k, k<j\vee k>j+1\}$. This results in a smaller problem of the form:
\begin{align}\label{eq:mals_projected_problem}
	\mathcal V_{j,2}^T \mathcal A_\text{TT} \mathcal V_{j,2} Y &= \mathcal V_{j,2}^T B_\text{TT}.
\end{align}
For s.p.d.\ operators $A_\text{TT}$, this much smaller problem typically has a better condition number
($\kappa\left(\mathcal V_{j,2}^T \mathcal A_\text{TT} \mathcal V_{j,2}\right) \le \kappa\left(\mathcal A_\text{TT}\right)$),
and its solution minimizes the error in the induced operator norm:
\begin{align*}
	Y &= \arg\min_Y \left\|\mathcal V_{j,2}Y - \mathcal X_\text{TT}^*\right\|_{\mathcal A_\text{TT}}^2 
	\qquad\text{for}\qquad Y= \left(\mathcal V_{j,2}^T\mathcal A_\text{TT}\mathcal V_{j,2}\right)^{-1}\mathcal V_{j,2}^T B_\text{TT},
\end{align*}
where $\mathcal X_\text{TT}^*$ denotes the true solution.
For non-symmetric operators $\mathcal A_\text{TT}$, this projection (Ritz-Galerkin) is often still successful \citep{Dolgov2014a}, but one might also consider a different projection for the left and the right-hand side of the operator (Petrov-Galerkin):
\begin{align*}
	\bar{\mathcal V}_j^T \mathcal A_\text{TT} \mathcal V_{j,2} Y &= \bar{\mathcal V}_j^T B_\text{TT}.
\end{align*}
Here, $\bar{\mathcal V}_j$ should have the same dimensions as $\mathcal V_{j,2}$ to ensure that projected problem is square (and thus usually easier to solve).
A possible non-symmetric approach is
\begin{align*}
	\bar{\mathcal V}_j Z &\approx \mathcal A_\text{TT} \mathcal V_{j,2} \qquad \text{with}\qquad \bar{\mathcal V}_j^T\bar{\mathcal V}_j = I.
\end{align*}
This is only possible approximately if $\bar{\mathcal V}_j$ should have low rank in the TT format again.
The solution of the projected problem then approximately minimizes the residual in the Frobenius norm (similar to GMRES).
In this paper, we will not further discuss this approach as we focus on the performance of the operations involved, but other variants of projections are possible.

\begin{algorithm}
\caption{TT-MALS}
\label{alg:TT-MALS}
\begin{algorithmic}[1]
\REQUIRE{Linear operator $\mathcal A_\text{TT}: \mathbf R^{n_1 \times \cdots \times n_d} \to \mathbf R^{n_1 \times \cdots \times n_d}$ and RHS $B_\text{TT} \in \mathbf R^{n_1 \times \cdots \times n_d}$,\\initial guess $X_\text{TT}\in\mathbf R^{n_1\times\cdots\times n_d}$, desired tolerance $\epsilon$, max. number of sweeps $m$}
\ENSURE{Approximate solution $X_\text{TT}$ with $\left\|B_\text{TT} - \mathcal A_\text{TT} X_\text{TT}\right\|_F \lesssim \epsilon \left\|B_\text{TT}\right\|_F$}
	\STATE {Right-orthogonalize $X_d \dots X_3$}
	\FOR {$i_\text{Sweep} = 1, \dots, m$}
		\STATE {$j_\text{start} \gets 1$ \textbf{if} $i_\text{Sweep} = 1$ \textbf{else} $2$}
		\FOR {$j = j_\text{start}, \dots, d-1$}
			\STATE {Left-orthogonalize $X_{j-1}$ \textbf{if} $j>1$}
			\STATE {Setup projection operator $\mathcal V_{j,2}$ using \eqref{eq:mals_projection_operator}}
			\STATE {Solve local problem $\mathcal V_{j,2}^T \mathcal A_\text{TT} \mathcal V_{j,2} Y = \mathcal V_{j,2}^T B_\text{TT}$ with initial guess $X_j\Join X_{j+1}$}
			\STATE {Update $X_j\Join X_{j+1} \gets Y$}
		\ENDFOR
		\IF {$\left\|B_\text{TT} - \mathcal A_\text{TT} X_\text{TT}\right\|_F \le \epsilon\|B_\text{TT}\|$}
			\RETURN
		\ENDIF
		\FOR {$j = d-2, \dots, 1$}
			\STATE {Right-orthogonalize $X_{j+2}$}
			\STATE {Setup projection operator $\mathcal V_{j,2}$ using \autoref{eq:mals_projection_operator}}
			\STATE {Solve local problem $\mathcal V_{j,2}^T \mathcal A_\text{TT} \mathcal V_{j,2} Y = \mathcal V^T B_\text{TT}$ with initial guess $X_j\Join X_{j+1}$}
			\STATE {Update $X_j\Join X_{j+1} \gets Y$}
		\ENDFOR
		\IF {$\left\|B_\text{TT} - \mathcal A_\text{TT} X_\text{TT}\right\|_F \le \epsilon \|B_\text{TT}\|$}
			\RETURN
		\ENDIF
	\ENDFOR
	\STATE {Abort: not converged in $m$ sweeps!}
\end{algorithmic}
\end{algorithm}
To solve the global problem, the MALS algorithm sweeps from the first two to the last two dimensions and back, see \autoref{alg:TT-MALS}.
This can be interpreted as a moving subspace correction as discussed in \cite{Oseledets2018}.
In contrast, the ALS algorithm only projects the problem onto the subspace of all but one sub-tensor.
This yields smaller local problems but provides no mechanism to increase the TT-rank of the approximate solution.
So in the simplest form, it can only converge for special initial guesses that already have the same TT-rank as the desired approximate solution.
In \autoref{sec:AMEn_method}, we discuss a possible way to avoid this problem.

\subsubsection{Inner solver: TT-GMRES}\label{subsubsec:inner_solver_tt_gmres}
The projected problem~\eqref{eq:mals_projected_problem} again has a structure similar to \eqref{eq:linear_system} but with just two dimensions. 
More specifically, after some contractions, the projected operator has the form
\begin{align*}
	\mathcal V_{j,2}^T \mathcal A_\text{TT} \mathcal V_{j,2} &= \bar A_{j-1,\text{left}}\Join A_j \Join A_{j+1} \Join \bar A_{j+2,\text{right}}
\end{align*}
with $\bar A_{j-1,\text{left}} \in \mathbf R^{r_{j-1}\times r_{j-1} \times r^A_{j-1}}$ and $\bar A_{j+2,\text{right}} \in\mathbf R^{r^A_{j+1}\times r_{j+1} \times r_{j+1}}$.
And the initial guess as well as the required solution is in the form:
\begin{align*}
	Y &= X_j \Join X_{j+1} \in\mathbf R^{r_{j-1} \times n_j \times n_{j+1} \times r_{j+1}}.
\end{align*}
Using tensor contractions, we can also express the projected right-hand side as:
\begin{align*}
	\mathcal V_{j,2}^T B_\text{TT} &= \bar B_{j,\text{left}} \Join \bar B_{j+1,\text{right}}.
\end{align*}
As long as the rank of $Y$ during the iteration is much smaller than $r_{j-1}n_j$, respectively $n_{j+1}r_{j+1}$, it is usually beneficial to use the factorized form for $Y$.
This results in an inner-outer scheme with an outer MALS and an inner TT-GMRES algorithm.
From our observation, this also yields slightly smaller ranks in the outer MALS than using standard GMRES with the dense form of $Y$ and a subsequent factorization.
\begin{remark}
	In the first MALS sweeps, it is not necessary to solve the inner problem very accurately.
	So one can use a larger relative tolerance for the inner iteration than for the outer iteration.
	The same yields in the last MALS sweeps (close to the solution).
	Combining both aspects, we employ a relative tolerance of:
	\begin{align*}
		\bar \epsilon_\text{inner} &= \max\left(\epsilon_\text{inner}, \epsilon\frac{\|B_\text{TT}\|}{\|\mathcal A_\text{TT} X_\text{TT} - B_\text{TT}\|}\right), \quad\text{with}\quad\epsilon_\text{inner} := \sqrt{\epsilon}.
	\end{align*}
\end{remark}
\begin{remark}
	One might wonder if an inner-outer iteration scheme with outer (flexible) TT-GMRES and inner MALS (as preconditioner) might also work.
	In our experiments, this results in super-linear (up to exponential) growth of the TT-ranks in the Arnoldi iteration.
	This can be explained by the fact that the displacement structure of the linear operator is not retained through this form of a varying preconditioner.
\end{remark}
For symmetric problems, our implementation switches to TT-MINRES and for positive definite problems, one can also employ a tensor-train variant of the CG algorithm.

\subsection{AMEn method}\label{sec:AMEn_method}
For the MALS method above, always two sub-tensors are considered at once in the inner problem.
One can also consider only one sub-tensor (ALS) at a time using the projection operator
\begin{align}\label{eq:als_projection_operator}
	\mathcal V_{j,1}&: \mathbf R^{r_{j-1}\times n_j \times r_j} \to \mathbf R^{n_1 \times \cdots \times n_d}, \quad \text{with} \nonumber\\
	\mathcal V_{j,1}Y&= X_1\Join \cdots \Join X_{j-1} \Join Y \Join X_{j+1} \Join \cdots \Join X_d.
\end{align}
This results in a smaller local operator. By contracting sub-tensors of $\mathcal V_{j,1}$ and $\mathcal A_\text{TT}$, one obtains:
\begin{align}\label{eq:amen_local_operator}
	\bar {\mathcal A}_\text{TT} &:= \mathcal V_{j,1}^T \mathcal A_\text{TT} \mathcal V_{j,1} = \bar A_{j-1,\text{left}} \Join A_j \Join \bar A_{j+1,\text{right}}.
\end{align}
However, one needs a way to adapt the rank of the approximate solution and to ensure convergence.
An early approach from physics was just to increase the rank through adding random directions.
A more sophisticated method is introduced in \cite{Dolgov2014a} named AMEn (alternating minimal energy).
The main idea is to enrich the subspace after each inner iteration with directions obtained from the current residual tensor.
\revised{For this, the current residual tensor is projected onto the subspace of the left sub-tensors of the current approximation when sweeping left-to-right, respectively the right sub-tensors when sweeping right-to-left.
The corresponding projection operators are given by:
\begin{align}\label{eq:amen_residual_projection_operator}
	\mathcal V_{j,d-j}&: \mathbf R^{r_{j-1}\times n_j\times\cdots\times n_d} \to \mathbf R^{n_1\times \cdots \times n_d}, \quad \text{with} &\mathcal V_{j,d-j}Y&= X_1 \Join \cdots \Join X_{j-1} \Join Y, \nonumber\\
	\mathcal V_{1,j}&: \mathbf R^{n_1\times\cdots\times n_j \times r_{j+1}} \to \mathbf R^{n_1\times \cdots \times n_d}, \quad \text{with} &\mathcal V_{1,j}Y&= Y\Join X_{j+1} \Join \cdots \Join X_d.
\end{align}
}%
For s.p.d.\ operators, this results in a steepest descent step for which global convergence is shown in \cite{Dolgov2014a}.
In practice, it is often sufficient to approximate a few directions of the residual tensor and add them to the current subspace to obtain fast convergence.
The standard form of the AMEn algorithm is depicted in \autoref{alg:TT-AMEn} (based on the SVD variant in \cite{Dolgov2014a}).
\begin{algorithm}[ht!]
\caption{TT-AMEn}
\label{alg:TT-AMEn}
\begin{algorithmic}[1]
\REQUIRE{Linear operator $\mathcal A_\text{TT}$, RHS $B_\text{TT}$ and initial guess $X_\text{TT}$, desired tolerance $\epsilon$,\\
         max. number of sweeps $m$, number of enrichment directions $k$}
\ENSURE{Approximate solution $X_\text{TT}$ with $\left\|B_\text{TT} - \mathcal A_\text{TT} X_\text{TT}\right\|_F \lesssim \epsilon \left\|B_\text{TT}\right\|_F$}
	\STATE {Right-orthogonalize $X_d \dots X_2$}
	\STATE {Calculate residual $R_\text{TT} \gets \mathcal A_\text{TT} X_\text{TT} - B_\text{TT}$}
	\STATE {Right-orthogonalize $R_d \dots R_2$ \revised{and let $Z_\text{TT} \gets R_\text{TT}$}}
	\FOR {$i_\text{Sweep} = 1, \dots, m$}
		\FOR {$j = 1, \dots, d-1$}
			\STATE {\revised{Setup projection operators $\mathcal V_{j,1}, \mathcal V_{j,d-j}$ using \eqref{eq:als_projection_operator} and \eqref{eq:amen_residual_projection_operator}}}
			\STATE {Solve local problem $\mathcal V_{j,1}^T \mathcal A_\text{TT} \mathcal V_{j,1} Y = \mathcal V_{j,1}^T B_\text{TT}$ with initial guess $X_j$}
			\STATE {Update $X_j\gets Y$}
			\STATE {Update $R_j$ s.t. $R_\text{TT} = \mathcal A_\text{TT} X_\text{TT} - B_\text{TT}$} \label{step:updateR}
			\IF {$\left\|R_\text{TT}\right\|_F \le \epsilon\|B_\text{TT}\|$}
				\RETURN
			\ENDIF
			\revised{
			\STATE {Update $Z_j$, s.t. $Z_j \Join \cdots \Join Z_d = \mathcal V_{j,d-j}^T R_\text{TT}$}\label{step:updateZ}
			\STATE {Left-orthogonalize $X_j$ (updating $X_{j+1}$ s.t. $X_j\Join X_{j+1}$ remains the same)}\label{step:OrthoSVD_Z}
			\STATE {Left-orthogonalize $Z_j$ w.r.t.\ $X_j$}\label{step:selectZ}
			}
			\STATE {Append $Z_{:,:,1:k} \Join 0$ to $X_j \Join X_{j+1}$}\label{step:appendZ}
		\ENDFOR
		\STATE {Similarly sweep from $d$ to $2$.}
	\ENDFOR
	\STATE {Abort: not converged in $m$ sweeps!}
\end{algorithmic}
\end{algorithm}
The step to update the residual \revised{$R_\text{TT}=\mathcal A_\text{TT} X_\text{TT}$} (line~\ref{step:updateR}) requires saving all intermediate matrices from left- respectively right-orthogonalization of the sub-tensors of $R_\text{TT}$.
To calculate $k$ suitable directions to enrich the subspace in the left-to-right sweep (line~\ref{step:selectZ}), we left-orthogonalize $Z_j$ w.r.t.\ $X_j$, such that:
\revised{
\begin{align*}\label{eq:amen_left_enrichment_directions}
	U{\color{gray}SV^T} &= (I - (X_j)_{\text{unfold}_l}(X_j)_{\text{unfold}_l}^T)(R_j)_{\text{unfold}_l}, \qquad  (Z_j)_{\text{unfold}_l} = U.
\end{align*}
This results in the method AMEn+SVD from \cite{Dolgov2014a}.
}%
To append the directions (line~\ref{step:appendZ}), we concatenate the tensors, such that $\bar r_j = r_j+k$ and $(\bar X_j)_{:,:,r_j+1:r_j+k} = Z_{:,:,1:k}$ and $(\bar X_{j+1})_{r_j+1:r_j+k,:,:} = 0$.
The same is done in the right-to-left sweep with mirrored dimensions.

As can be seen from \autoref{alg:TT-AMEn}, calculating the required directions from the residual mostly involves updating the sub-tensors of the residual in every step of the sweep.
This is not significantly more work than calculating the residual in the first place.
Still, the complete algorithm is so cheap that the residual calculation accounts for a significant part of the overall runtime.
That is why \cite{Dolgov2014a} discusses several more heuristic ways to determine suitable enrichment directions.

\subsubsection{\revised{TT-AMEn+ALS}}
\revised{
The most promising variant from \cite{Dolgov2014a} is based on an ALS-like approximation of the residual.
As the resulting complete algorithm is not explicitly shown there, we illustrate the required steps \autoref{alg:TT-AMEn+ALS} and discuss important implementation details.
}

First, to really decrease the work, one needs a way to check the error without using the global residual norm.
In \cite{Dolgov2014a}, this is not further discussed but in the code used for the numerical experiments of \cite{Dolgov2014a}, the global residual error is estimated using the projected residuals (before solving the local problem), assuming that
\begin{align*}
	\frac{1}{2\sqrt{d-1}} \max_j\left(\left\|\mathcal V_{j,1}^T(\mathcal A_\text{TT} X_\text{TT} - B_\text{TT})\right\|_F)\right) &\lesssim \left\|\mathcal A_\text{TT} X_\text{TT} - B_\text{TT}\right\|.
\end{align*}
The factor $1/2$ is just a heuristic way to ensure that the global residual norm is smaller than the tolerance.

\revised{
Second, one needs a cheap way to enrich the subspace.
For this, a rough approximation of the residual is usually sufficient which can be obtained by a fixed-rank ALS iteration ($\operatorname{ALS}(l)$).
This results in the following update of the $j$-th subtensor of the current approximation of the residual $\tilde R_\text{TT}$:
}%
\begin{align}
	\tilde R_j' &:= \mathcal W_{j,1}^T\left(\mathcal A_\text{TT} X_\text{TT} - B_\text{TT}\right) = \mathcal W_{j,1}^T \mathcal A_\text{TT} X_\text{TT} - \mathcal W_{j,1}^T B_\text{TT}, \quad\text{where}\\
	\mathcal W_{j,1} &:\mathbf R^{l\times n_j \times l} \to \mathbf R^{n_1\times\cdots\times n_d}, \quad \mathcal W_{j,1}Y=\tilde R_1\Join \cdots \tilde R_{j-1} \Join Y \Join \tilde R_{j+1} \Join \cdots \Join \tilde R_d. \label{eq:amen_als_projection_operator}
\end{align}
\revised{
Here, $\tilde R_1, \dots, \tilde R_{j-1}$ must be left-orthogonal and $\tilde R_{j+1},\dots,\tilde R_d$ right-orthogonal.
We remark that the approximation $\tilde R_\text{TT}$ cannot be used to check convergence as $\|\tilde R_\text{TT}\|_F\le \|R_\text{TT}\|_F$.
As for the projected local problem, the required additional terms can be updated successively in each step of the sweep.
}
\begin{algorithm}
\caption{\revised{TT-AMEn+ALS}}
\label{alg:TT-AMEn+ALS}
\begin{algorithmic}[1]
\REQUIRE{Linear operator $\mathcal A_\text{TT}$, RHS $B_\text{TT}$, initial guess $X_\text{TT}$, desired tolerance $\epsilon$,\\
         max. number of sweeps $m$, enrichment rank $k$, inner tolerance $\epsilon_\text{inner}$, \revised{approx. residual rank $l$}}
\ENSURE{Approximate solution $X_\text{TT}$ with $\left\|B_\text{TT} - \mathcal A_\text{TT} X_\text{TT}\right\|_F \lesssim \epsilon \left\|B_\text{TT}\right\|_F$}
	\STATE {Right-orthogonalize $X_d \dots X_2$}
	\STATE {\revised{Initialize approx. residual $\tilde R_1, \dots, \tilde R_d$ with rank $l$}}
	\STATE {\revised{Right-orthogonalize $\tilde R_d \dots \tilde R_2$ and let $\tilde Z_\text{TT} \gets \tilde R_\text{TT}$}}
	\FOR {$i_\text{Sweep} = 1, \dots, m$}
		\FOR {$j = 1, \dots, d$}
			\STATE {\revised{Setup projection operators $\mathcal V_{j,1}, \mathcal V_{j,d-j}$ using \eqref{eq:als_projection_operator} and \eqref{eq:amen_residual_projection_operator}}}
			\IF {$i_\text{Sweep} = 1 \vee j \ne 1$}
				\STATE {Calculate local residual norm: $\delta_j \gets \| \mathcal V_j^T(\mathcal A_\text{TT}X_\text{TT}-B_\text{TT})\|_F$}
				\STATE {Adapt tolerance: $\bar\epsilon_\text{inner,abs} \gets \max(\delta_j\,\epsilon_\text{inner},\; \|\mathcal V_j^T B_\text{TT}\|_F\,\epsilon/(2\sqrt{d-1}) )$}
				\STATE {Approximately solve local problem $\mathcal V_{j,1}^T \mathcal A_\text{TT} \mathcal V_{j,1} Y = \mathcal V_{j,1}^T B_\text{TT}$\\
						with initial guess $X_j$ up to abs. tolerance $\bar\epsilon_\text{inner,abs}$}
				\STATE {Update $X_j\gets Y$}
				\STATE {\revised{Setup residual projection operator $\mathcal W_{j,1}$ using \eqref{eq:amen_als_projection_operator}}}
				\STATE {\revised{Update $\tilde R_j\gets \mathcal W_{j,1}^T(\mathcal A_\text{TT}X_\text{TT}-B_\text{TT})$}}
			\ENDIF
			\IF {$j < d$}
				\STATE {\revised{Update $\tilde Z_j$, s.t. $\tilde Z_j \Join \cdots \Join \tilde Z_d = \mathcal V_{j,d-j}^T \tilde R_\text{TT}$}}
				\STATE {\revised{Left-orthogonalize $\tilde R_j$ (updating $\tilde R_{j+1}$ s.t. $\tilde R_j\Join \tilde R_{j+1}$ remains the same)}}
				\STATE {Left-orthogonalize $X_j$ \revised{(updating $X_{j+1}$ s.t. $X_j\Join X_{j+1}$ remains the same)}}
				\STATE {Left-orthogonalize $\tilde Z_j$ w.r.t.\ $X_j$}
				\STATE {Append $\tilde Z_{:,:,1:k} \Join 0$ to $X_j \Join X_{j+1}$}
			\ENDIF
		\ENDFOR
		\IF {$\max_j(\delta_j) \le \frac{\epsilon}{2\sqrt{d-1}}$}
			\RETURN
		\ENDIF
		\STATE {Similarly sweep from $d$ to $1$ and check convergence.}
	\ENDFOR
	\STATE {Abort: not converged in $m$ sweeps!}
\end{algorithmic}
\end{algorithm}
\revised{
In our implementation, we use $B_\text{TT}$ as initial guess for the approximate residual $\tilde R_\text{TT}$ and truncate it, respectively extend it by random directions to obtain the desired rank $l$.
}%
In \autoref{subsec:performance_of_tt_amen}, we show numerical experiments with both the full AMEn and the AMEn+ALS algorithm.

\subsection{Preconditioning}
For iterative solvers of linear systems, it is a common approach to employ a preconditioner to obtain much faster convergence.
Of course, we can precondition all previously discussed algorithms.
However, some additional aspects should be considered when preconditioning linear solvers in the TT format.
On the one hand, for MALS and AMEn, we can employ \emph{local preconditioners} for the projected operators $\mathcal V_j^T \mathcal A_\text{TT} \mathcal V_j$.
In this case,
\begin{itemize}
	\item one needs to calculate a different preconditioner in every step of the sweep,
	\item often only a few local iterations are performed,
	\item one may need to contract $\mathcal V_j^T \mathcal A_\text{TT} \mathcal V_j = \bar A_{j-1,\text{left}} \Join A_j \Join \bar A_{j+1,\text{right}}$ (costly).
\end{itemize}
On the other hand, we could employ a \emph{global preconditioner} which is
\begin{itemize}
	\item either applied to tensor-train ``vectors'' (TT-GMRES),
	\item or directly to the tensor-train operator $\mathcal A_\text{TT}$ (TT-GMRES, MALS, AMEn),
	\item but is not tailored to the local problems for MALS/AMEn.
\end{itemize}
Furthermore, a \emph{global preconditioner} (and also \emph{local precondioners} for MALS) should retain the low-rank of the solution and right-hand side.
To emphasize this: just the fact that $\kappa(\mathcal P \mathcal A_\text{TT})\ll\kappa(\mathcal A_\text{TT})$ does not imply that $\operatorname{rank}_\text{TT}(\mathcal PB_\text{TT}) \le \operatorname{rank}_\text{TT}(B_\text{TT})$ for left-preconditioning with $\mathcal P$.
And similarly the same problem occurs for the approximate solution $X_\text{TT}$ for right-preconditioning.

\subsubsection{TT-rank-1 preconditioner}
Considering the desired properties, we suggest a simple, \emph{global} rank-1 preconditioner that usually is cheap to calculate (for operators of sufficiently low rank).
The two-sided variant can be constructed as follows:
First, approximate the operator with a tensor-train of rank~1 using TT-truncation: $\tilde{\mathcal A}_\text{TT} = \tilde A_1 \otimes \tilde A_2 \otimes \dots \otimes \tilde A_d = \operatorname{trunc}_{r=1}(\mathcal A_\text{TT})$.
Then, perform an SVD for each sub-matrix ($U_j S_j V_j^T = \tilde A_j$, $j = 1, \dots, d$) to construct the TT operators for the left and the right:
\begin{align*}
	\mathcal P_\text{TT,left} &= (S_1^{-1/2}U_1^T) \otimes \cdots \otimes (S_d^{-1/2}U_d^T),&
	\mathcal P_\text{TT,right} &= (V_1 S_1^{-1/2}) \otimes \cdots \otimes (V_d S_d^{-1/2}).
\end{align*}
This yields the preconditioned system:
\begin{align*}
	\left(\mathcal P_\text{TT,left} \mathcal A_\text{TT} \mathcal P_\text{TT,right}\right) Y_\text{TT} &= \mathcal P_\text{TT,left} B_\text{TT} \qquad\text{with}\qquad X_\text{TT} = \mathcal P_\text{TT,right} Y_\text{TT}.
\end{align*}
For a symmetric operator $\mathcal A_\text{TT}$ where for each sub-tensor $(A_j)_{:,k,l,:}=(A_j)_{:,l,k,:}$, the preconditioned operator is still symmetric.
And as the preconditioner has rank one, it does not lead to higher ranks for the right-hand side or the exact solution $X_\text{TT}$.
However, if the operator has a displacement structure, this is not preserved for the preconditioned operator.

\section{Comparison of algorithms}\label{sec:comparison_of_methods}
In the following we present numerical experiments for linear systems in TT format.
We use the pitts library \citep{RoehrigZoellner2024pitts} which also contains \revised{the setup and output for all results shown in this paper.}
\revised{For illustrating the numerical behavior, we consider} a simple multidimensional convection-diffusion equation.
It is discretized using a finite difference stencil, e.g., in 1d the operator is
\begin{align*}
	\frac{\operatorname{tridiag}(-1,2,-1)}{h_j^2} + \frac{c}{\sqrt{d}}\;\frac{\operatorname{tridiag}(0,1,-1)}{h_j}, \qquad\text{for } h_j=\frac{1}{n_j+1}.
\end{align*}
Here, $c$ denotes a convection constant and the convection direction is diagonal through all dimensions.
As right-hand side, we use a vector of all ones (rank~1) or just a tensor-train with chosen rank and random sub-tensors.
All cases use double-precision calculations and a desired residual tolerance of $\epsilon=10^{-8}$.

\subsection{Behavior of tensor-train ranks in the calculation}\label{subsec:behavior_of_tt_ranks}
First, to understand the computational complexity of the different methods, we show the behavior of the TT ranks during the calculation.
\begin{figure}
\subcaptionbox{Krylov basis ranks for TT-GMRES\label{fig:tt-gmres_ranks}}{\includegraphics{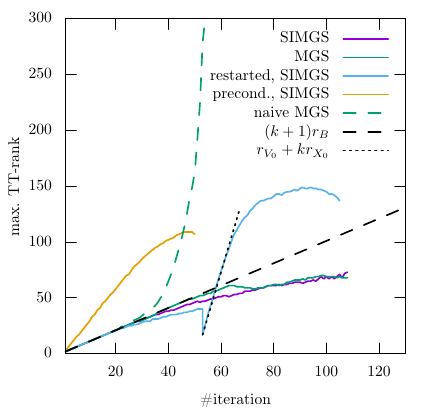}}
\hfill
\subcaptionbox{Ranks for precond.\ MALS with inner TT-GMRES\label{fig:mals_ranks}}{\includegraphics{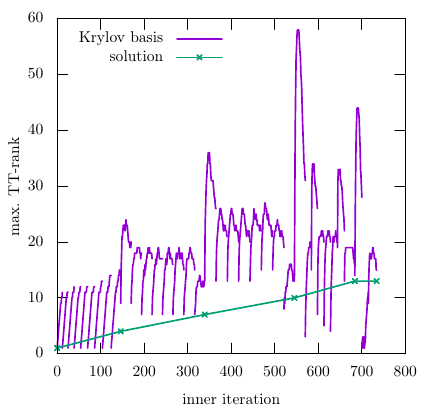}}
\caption{Tensor-train ranks for the Krylov basis, respectively the approximate solution for a $20^{10}$ convection-diffusion\ problem ($c=10$) and RHS $B_\text{TT}$ of ones.
\revised{%
For TT-GMRES (left), both MGS variants lead to inaccurate solutions that are not within the desired residual tolerance in contrast to all cases with SIMGS.
Overall, more accurate orthogonalization (SIMGS) without restart and preconditioning features the lowest maximal ranks during the calculation.
For MALS (right), the solution ranks only increase slowly with each sweep (as intended), but the Krylov basis vectors of the inner iteration again yield higher ranks.
}}
\end{figure}
As shown in~\autoref{fig:tt-gmres_ranks}, the improved orthogonalization in the TT-GMRES algorithms (SIMGS, \autoref{alg:SIMGS}) results in slightly slower rank growth in the Krylov basis than standard MGS; Both show the expected linear growth (for operators with a displacement structure).
A naive MGS implementation (all truncations performed with tolerance $\delta_i$ in \autoref{alg:TT-GMRES}) results in exploding ranks after a few iterations.
\revised{%
Furthermore, the calculated solutions of the MGS variants are not within the desired residual tolerance (with $\|\mathcal A_\text{TT}X_\text{TT}-B_\text{TT}\|_F \approx 10 \epsilon$).
}
With the TT-rank-1 preconditioner, the algorithm converges after less than half the number of iterations, but higher ranks occur (as the preconditioned operator does not have a displacement structure).
Similarly, a GMRES restart also results in higher TT ranks.
Overall, the TT-GMRES has a computational complexity that is cubic in the \revised{maximal} TT rank \revised{($r_\text{max}$)} of the Krylov basis: $O(dnr_\text{max}^3 + dn^2r_\text{max}^2(r^A)^2)$.
For the MALS method, we observe a (linearly) growing rank of the approximation solution in the iteration up to the rank \revised{($r$)} of the resulting approximate solution, see \autoref{fig:mals_ranks}.
However, in the inner TT-GMRES iterations, the Krylov basis ranks behave like after a GMRES restart after the first sweep.
So significantly higher ranks occur in the inner iteration.
The effect on the computational complexity is only quadratic for MALS (only two sub-tensors in the inner problem): $O(dnrr_\text{max}^2 + dn^2rr_\text{max}(r^A)^2)$.
The behavior of the rank of the approximate solution for AMEn is similar to MALS.
But as the local problem consists only of one sub-tensor, no higher intermediate ranks can occur \revised{in the inner iteration}.
So only $O(dnr^3 + dn^2r^2(r^A)^2)$ operations are needed \revised{assuming that the enrichment rank $k$ is chosen appropriately.}

\subsection{Computational complexity of the different methods}
\begin{figure}
\subcaptionbox{Varying dimensions (RHS $B_\text{TT}$ of ones)\label{fig:tt_method_comparison:dims}}{\includegraphics{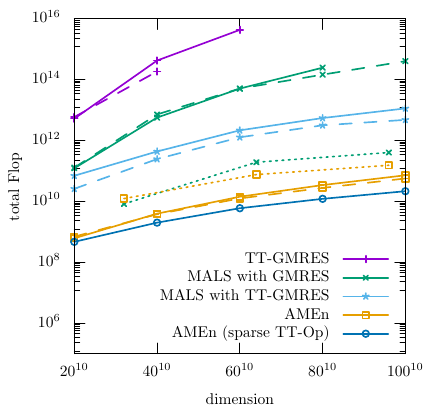}}
\hfill
\subcaptionbox{Random RHS $B_\text{TT}$ with varying ranks (dimension $20^{10}$)\label{fig:tt_method_comparison:rhs}}{\includegraphics{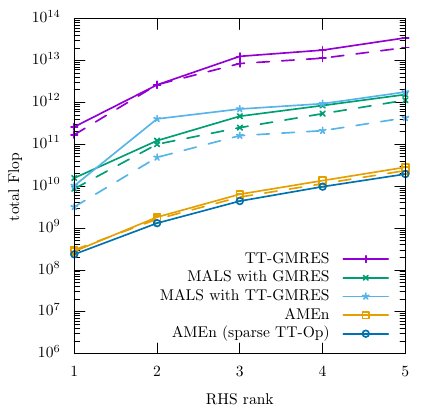}}
\caption{Number of floating-point operations measured using likwid \citep{Treibig2010} for a convection-diffusion problem ($c=10$). Dashed lines use the TT-rank-1 preconditioner. Dotted lines first transform the problem to the QTT format \citep{Khoromskij2011}.
\revised{%
In all cases, AMEn requires orders of magnitude fewer operations than MALS and TT-GMRES.
}}\label{fig:tt_method_comparison}
\end{figure}

The TT ranks explain the results in \autoref{fig:tt_method_comparison}: TT-GMRES needs 10-100$\times$ more operations than MALS. MALS needs 10-100$\times$ more operations than AMEn.
For larger individual dimensions ($n_j$), MALS with inner TT-GMRES needs fewer operations than with inner MALS.
In addition, applying the operator becomes more costly, and it is beneficial to use a sparse matrix format for the sub-tensors of the operator (sparse TT-Op variant in \autoref{fig:tt_method_comparison:dims}).
We also measured the behavior using the quantics TT format (QTT, \cite{Khoromskij2011}) where we just convert the operator to a $2^{\bar d}$ tensor.
For our test case here, QTT-MALS needs fewer operations than TT-MALS (ranks in inner iteration can grow at most by a factor of~2).
However, QTT-AMEn needs more operations than TT-AMEn here.
Looking at the approximate solution, it has TT ranks $(\dots,r,2r,3r,4r,\dots,4r,3r,2r,r,\dots)$ where $r$ is the rank of the approximate solution without transformation to QTT format.
So the QTT format is not beneficial here as the solution is not well approximable with small QTT ranks.

\section{Performance of algorithmic building blocks}\label{sec:building_blocks}
In the following, we discuss the required basic tensor-train operations (``building blocks'') for the algorithms from \autoref{sec:numerical_algorithms}.
We focus on the node-level performance on a single multicore CPU with some remarks for a distributed parallel implementation.
We only consider those operations where we see a significant improvement over the standard implementation as introduced in \cite{Oseledets2011}.
These are in particular left-/right-orthogonalization with/without truncation, TT addition with subsequent truncation (TT-axpby+trunc), and faster contractions. 

\subsection{Replacing costly SVDs and pivoted QR decompositions with faster but less accurate alternatives}\label{sec:building_blocks_svd_and_qr}
In \cite{RoehrigZoellner2022}, the authors present a significantly faster implementation for decomposing a large dense tensor in the TT format (TT-SVD) using a Q-less tall-skinny QR (TSQR) algorithm.
For the TT-SVD, one starts with a decomposition of the form:
\begin{align*}
	&&&\min_{B,Q^T}\left\|(X)_{\text{unfold}_{d-1}} - BQ^T \right\|_F && \text{with}& Q^TQ&=I.
\end{align*}
$(X)_{\text{unfold}_{d-1}}$ is a tall-skinny matrix, and the standard algorithm uses a truncated SVD to build $B$ and $Q$:
\begin{align*}
	&&(X)_{\text{unfold}_{d-1}} &= U S V^T && \text{and}& B&=U_{:,1:r_1} S_{1:r_1,1:r_1}, \quad Q=V_{:,1:r_1}.
\end{align*}
The optimized algorithm uses the following steps instead (grayed-out matrices are not calculated):
\begin{align*}
	(X)_{\text{unfold}_{d-1}} &= {\color{gray} \bar Q} R, \quad R \approx {\color{gray}\bar U} S V^T, &&\text{and}& B&=(X)_{\text{unfold}_{d-1}}V_{:,1:r_1}, \quad Q=V_{:,1:r_1} .
\end{align*}
As $V$ has orthonormal columns, the matrix $B$ can be calculated accurately this way.

For the left-/right-orthogonalization (see first loop in TT-rounding in \cite{Oseledets2011}), we unfortunately need slightly different operations, e.g., in the left-to-right QR sweep:
\begin{align*}
	(X_j)_\text{left} &= QB, &&\text{and}& (X_j')_\text{left} &=Q, \quad X_{j+1}' = B \Join X_{j+1}.
\end{align*}
And very similarly for a left-to-right SVD sweep:
\begin{align*}
	&\min_{Q,B}\|(X_j)_\text{left} - QB\|_F &&\text{and}& (X_j')_\text{left} &= Q, \quad X_{j+1}' = B \Join X_{j+1}.
\end{align*}
The only difference is that the pure orthogonalization is exact up to the numerical rank whereas the truncation intentionally cuts off with a given tolerance.
We can use the same trick as for the TT-SVD in a slightly different way here.
For the orthogonalization, one obtains with pivoting matrix $P$:
\begin{align}\label{eq:inaccurate_orthogonalization}
	(X_j)_\text{left}P &= {\color{gray} Q}R && \text{and}& \tilde X_j' &= X_j \Join (PR^{-1}), \quad X_{j+1}' = (RP^T) \Join X_{j+1} .
\end{align}
And again similarly for the truncation:
\begin{align}\label{eq:inaccruate_truncation}
	(X_j)_\text{left} &= {\color{gray} Q}R, \quad R= {\color{gray} U}SV^T && \text{and}& \tilde X_j' &= X_j \Join (VS^{-1}), \quad X_{j+1}' = (SV^T) \Join X_{j+1} .
\end{align}
There is a difference to the TT-SVD here:
in both cases, the faster formulation might introduce a significant numerical error ($\|\tilde X_j' - X_j'\|_F$). For $\kappa(R)\gg 1$, respectively $\kappa(S)\gg 1$,
\revised{%
applying $R^{-1}$, respectively $VS^{-1}$, can be numerically unstable.
}%
For the SVD-sweep, one can see that the error is only large in ``unimportant'' directions and should not affect the accuracy of the complete tensor-train:
Assuming $X_1,\dots,X_{j-1}$ are left-orthogonal and $X_{j+1},\dots,X_d$ right-orthogonal, then:
\begin{align*}
	\|\tilde X_\text{TT} - X_\text{TT} \|_F &= \| (\tilde X_j' - X_j') \Join X_{j+1}'\|_F = \| (\tilde X_j' - X_j') \Join S \|_F.
\end{align*}
This is not ensured for the orthogonalization step. 
However, this shows that the error should be weighted by the singular values of the corresponding unfolding of the approximated tensor.
So for the truncation, the left-to-right QR sweep is followed by a right-to-left SVD sweep with:
\begin{align*}
	\tilde X_j'' &= \tilde X_j' \Join (\bar U \bar S), &&\text{for}& \revised{\bar U \bar S \bar V^T} &\revised{= (\bar X_{j+1})_\text{right}.}
\end{align*}
We suggest checking the error a posteriori using the orthogonalization error weighted by the singular values of the $j$-th unfolding:
\begin{align*}
	\left\|\bar S^2 - (\tilde X_j'')_\text{left}^T (\tilde X_j'')_\text{left}\right\|_\infty &\overset{?}{\approx} \epsilon .
\end{align*}
In our numerical tests this was always the case even if $\|I - (\tilde X_j')_\text{left}^T (\tilde X_j')_\text{left}\|_\infty \gg \epsilon$.
If this a posteriori check fails, one can still recalculate the orthogonalization with the standard algorithm.

\subsection{Exploiting orthogonalities in TT-axpby+trunc}\label{sec:tt_axpby_exploiting_orthogonalities}
To add tensors in TT format ($Z_\text{TT} = \alpha X_\text{TT} + \beta Y_\text{TT}=Z_1\Join\cdots\Join Z_d$), one combines their sub-tensors:
\begin{align*}
	(Z_1)_{1,:,:} &= \begin{pmatrix} (X_1)_{1,:,:} & (Y_1)_{1,:,:}\end{pmatrix}, \\
	(Z_j)_{:,i,:} &= \begin{pmatrix} (X_j)_{:,i,:} & 0 \\ 0 & (Y_j)_{:,i,:} \end{pmatrix}, \quad\text{for}\quad i=1,\dots,n_j, \; j=2,\dots,d-1 \\
	(Z_d)_{:,:,1} &= \begin{pmatrix} \alpha (X_d)_{:,:,1} \\ \beta (Y_d)_{:,:,1} \end{pmatrix} .
\end{align*}
This operation is needed in the TT-GMRES algorithm as well as in MALS and the standard AMEn.
For AMEn, the operation is performed step-by-step in each step of a sweep.
In all cases, the sub-tensors of $X_\text{TT}$ and $Y_\text{TT}$ are already left- or right-orthogonal.
And after the addition, one performs a left- or right-orthogonalization of the $Z_\text{TT}$ for a subsequent truncation step.
We can exploit the block non-zero structure and pre-existing orthogonalities.
In the following, we assume $\operatorname{rank}_\text{TT}(X_\text{TT}) \ge \operatorname{rank}_\text{TT}(Y_\text{TT})$ and that $X_1,\dots,X_{d-1}$ are left-orthogonal.
Then, we can compute left-orthogonal sub-tensors for $\bar Z_1, \dots, \bar Z_{d-1}$ with smaller QR-decompositions $Q_1R_1$, \dots, $Q_{d-1}R_{d-1}$ than for the standard algorithm.
For the first sub-tensor, we obtain:
\begin{align}\label{eq:tt_axpby_ortho_1}
	&&(Z_1)_{1,:,:} &= \underbrace{\begin{pmatrix} (X_1)_{1,:,:} & Q_1 \end{pmatrix}}_{(\bar Z_1)_{1,:,:}} \begin{pmatrix}I & M_1 \\ 0 & R_1 \end{pmatrix}&&\text{with}& M_1 &= (X_1)_\text{left}^T (Y_1)_\text{left}, \nonumber\\[-1.2\baselineskip]
	&&              &                                                                                                   &&           & Q_1 R_1 &= (I-(X_1)_\text{left}(X_1)_\text{left}^T)(Y_1)_\text{left} .
\end{align}
For the next sub-tensors $j=2,\dots,d-1$, we obtain:
\begin{align}\label{eq:tt_axpby_ortho_middle}
	&& (Z_j')_{:,i,:} &= \begin{pmatrix} I & M_{j-1} \\ 0 & R_{j-1} \end{pmatrix} \begin{pmatrix} (X_j)_{:,i,:} & 0 \\ 0 & (Y_j)_{:,i,:} \end{pmatrix}&&\text{for}& i&=1,\dots,n_j\nonumber\\
	&&                &= \underbrace{\begin{pmatrix} \begin{matrix}(X_j)_{:,i,:} \\ 0\end{matrix} & (Q_j)_{:,i,:} \end{pmatrix}}_{(\bar Z_j)_{:,i,:}} \begin{pmatrix} I & M_j \\ 0 & R_j \end{pmatrix}&&\text{with}&
	                               M_j &= \bar X_j^T \bar Y_j, \\[-1.5\baselineskip]
	&&				  & &&&    Q_j R_j &= (I - \bar X_j \bar X_j^T) \bar Y_j, \nonumber
\end{align}
where we simplified the notation by introducing:
\begin{align*}
	&& \bar X_j &:= \left( \begin{pmatrix} I \\ 0 \end{pmatrix} \Join X_j \right)_\text{left}, & \bar Y_j &= \left( \begin{pmatrix} M_{j-1} \\ R_{j-1} \end{pmatrix} \Join Y_j \right)_\text{left} .
\end{align*}
The last sub-tensor simply results in:
\begin{align}\label{eq:tt_axpby_ortho_d}
	(\bar Z_d)_{:,:,1} &= \begin{pmatrix} I & M_{d-1} \\ 0 & R_{d-1} \end{pmatrix} \begin{pmatrix} \alpha (X_d)_{:,:,1} \\ \beta (Y_d)_{:,:,1} \end{pmatrix}.
\end{align}

\subsubsection{Stable residual calculation with inaccurate orthogonalization}
For (standard) AMEn and MALS, we update a left- respectively right-orthogonal representation of $\mathcal A_\text{TT} X_\text{TT}=:Y_\text{TT}$ in each step and calculate the residual $B_\text{TT} - Y_\text{TT}$ from it reusing the orthogonality.
This is susceptible to numerical errors as for the solution all directions should cancel out.
\revised{%
More specifically, we observed relative errors of up to $\epsilon={\sim}10^{-3}$ when the residual norm was close to the machine precision.
}%

\revised{%
In general, when combining inaccurate orthogonalization with the optimized TT-axpby+trunc algorithm, we need to adjust the formulas in \eqref{eq:tt_axpby_ortho_1} and \eqref{eq:tt_axpby_ortho_middle} to compensate for the loss of orthogonality.
Assuming $\bar X_j^T \bar X_j\approx I$, we suggest to employ a Cholesky decomposition to correct the projection in \eqref{eq:tt_axpby_ortho_middle} (and similarly in \eqref{eq:tt_axpby_ortho_1}) by using:
}%
\begin{align}
\revised{Q_j' R_j' }&= \revised{\left(I - \bar X_j (L_jL_j^T)^{-1} \bar X_j^T\right)\bar Y_j, \quad\text{with}\quad L_jL_j^T=\bar X_j^T \bar X_j.}
\end{align}
\revised{%
This emulates re-orthogonalization but is faster if $\bar Y_j$ has fewer columns than $\bar X_j$ as re-orthogonalization would require to update the matrix $\bar X_j'=\bar X_j(L_j)^{-T}$ instead.
}

\subsection{Faster contractions: inner iteration of AMEn}\label{sec:faster_contractions}
\revised{In AMEn, one needs to apply the projected TT operator from \eqref{eq:amen_local_operator} to a dense tensor:}
\begin{align*}
	Z &= \bar {\mathcal A}_{\text{TT},j} Y&&\text{with}& \bar {\mathcal A}_{\text{TT},j} = \bar A_{j-1,\text{left}} \Join A_j \Join \bar A_{j+1,\text{right}}, \qquad Y,Z \in \mathbf R^{r_{j-1}\times n_j \times r_j} .
\end{align*}
We can calculate this with the following contractions assuming that the ordering of dimensions denoted by $:$ is equal on the left- and right-hand sides of the equations:
\begin{align*}
	(Y')_{:,:,:,:} &\gets \sum_i (\bar A_{j+1,\text{right}})_{:,:,i} Y_{:,:,i},\\
	(Y'')_{:,:,:,:}&\gets \sum_{i_1,i_2}(A_j)_{:,:,i_1,i_2} (Y')_{i_2,:,:,i_1},\\
	Z_{:,:,:} &\gets \sum_{i_1,i_2} (\bar A_{j-1,\text{left}})_{:,i_1,i_2}(Y'')_{i_2,:,:,i_1} .
\end{align*}
If we reorder the array dimensions of the sub-tensors of $\bar {\mathcal A}_\text{TT}$ appropriately \revised{and let $(\hat A_1, \hat A_2, \hat A_3)$ denote the reordered tensors, we can use the following steps instead:}
\begin{align*}
	&&(\hat Y')_{:,:,:,k} &\gets \sum_i (\hat A_3)_{:,i,k} Y_{:,:,i} &&\text{for}& k&=1,\dots,r^A_j , \\
	&&(\hat Y'')_{:,:,:,k} &\gets \sum_{i_1,i_2} (\hat A_2)_{:,i_1,i_2,k} (\hat Y')_{:,:,i_1,i_2} &&\text{for}& k&=1,\dots,r^A_{j-1} , \\
	&&Z_{:,:,:} &\gets \sum_{i_1,i_2} (\hat A_1)_{:,i_1,i_2}(\hat Y'')_{:,:,i_1,i_2} .
\end{align*}
The benefit of this reordering is that we can combine the contracted and ``free'' dimensions ($(i_1,i_2)$ and $(:,:)$) to obtain fewer large dimensions which reduces overhead in the implementation.
In addition, we employ a column-major storage with a padding to obtain array strides that are multiples of the cache line length but not high powers of~2 to avoid cache thrashing.
Another approach consists in using optimized tensor contractions as discussed in detail in \cite{Springer2018} which may reduce the overhead due to several small dimensions.

\subsection{Resulting building block performance}
\revised{Here, we show experiments} on a single CPU socket of an AMD EPYC 7773X (``Zen 3 V-Cache'') with 64~cores and the Intel~oneMKL \citep{IntelOneMKL} as underlying BLAS / LAPACK library.\footnote{We also obtain qualitatively similar results on a 16-core Intel Xeon Gold 6246. We use a workaround for running Intel~oneMKL on AMD hardware defining the function \textsf{mkl\_serv\_intel\_cpu\_true} as discussed in \url{https://danieldk.eu/Posts/2020-08-31-MKL-Zen.html}.} 
We use the Q-less TSQR implementation discussed in \cite{RoehrigZoellner2022}.
\begin{figure}
\subcaptionbox{Timings for the truncated addition of two tensor-trains (dim.\ $50^{10}$, ranks $r_X=50$ and $r_Y=1,\dots,700$).\label{fig:axpby}}{\includegraphics{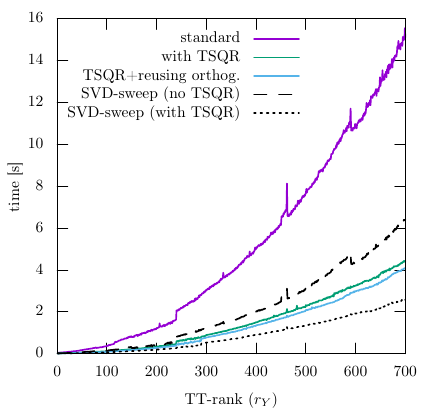}}
\hfill
\subcaptionbox{Performance of the contractions for multiplying a 3d TT-operator with a dense tensor (dim.: $r\times 50 \times r$)\label{fig:apply_op_dense}}{\includegraphics{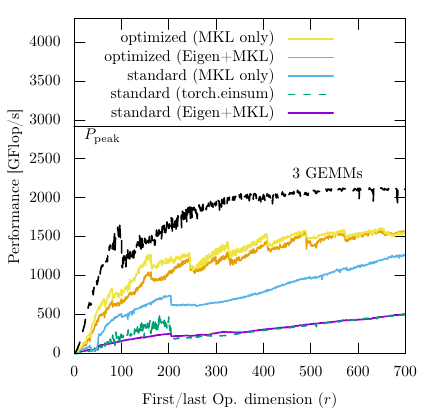}}
\caption{Effect of building block optimizations: 
\revised{%
For adding two tensors in the tensor-train format (left), we obtain a speedup of ${\sim}3.5$ by mapping the calculation onto faster linear algebra operations as explained in  \autoref{sec:building_blocks_svd_and_qr} and \autoref{sec:tt_axpby_exploiting_orthogonalities}.
For applying the linear operator of the inner problem in AMEn (right), we obtain a speedup of ${\sim}3$ through directly calling optimized BLAS routines and through reordering array dimensions.
}}
\end{figure}

\autoref{fig:axpby} shows timings for different variants of the TT-axpby+trunc operation.
We observe significant speedup through replacing SVDs and pivoted QR decompositions by a fast Q-less TSQR implementation:
\revised{%
For the QR orthogonalization sweep alone, the runtime improves by a factor of ${\sim}4.5$ for the largest case (difference between colored an black lines).
For the subsequent SVD sweep, the runtime improves by a factor of ${\sim}2.5$.
Overall speedup is about ${\sim}3.5$.
}%
We only see a small effect through reusing the orthogonality in the example here.
It uses random tensor-trains as inputs, so the rank of the result is the sum of the individual ranks.
In practice, e.g., for the Arnoldi iteration or for calculating the residual in AMEn, the rank of the resulting tensor-train is often smaller which leads to less work in the SVD sweep and thus a bigger effect through reusing orthogonalities in the preceding QR sweep.

For the contraction in the inner iteration of AMEn, we illustrate the performance of different variants in \autoref{fig:apply_op_dense}.
The operation consists of 3 tensor contractions, one of them is memory-bound and the other two are compute-bound for the chosen dimensions (for $r \gtrsim 100$).
The dotted line shows the performance of 3 equivalent matrix-matrix multiplications (GEMM).
In particular for small dimension $r$, there is a significant improvement through reordering and combining the dimensions as discussed in \autoref{sec:faster_contractions}.
Without appropriate padding, there are performance drops whenever the $r$ is close to a number dividable by, e.g., 128 (not shown here).
Overall, there is also a significant difference between using Eigen \citep{Eigen} with MKL as backend and directly calling MKL through the cblas interface, possibly because Eigen explicitly initializes the result to zero.
The standard implementation used here loops over a lot of possibly small GEMM operations.
\revised{We also show results using the function \texttt{torch.einsum} in pytorch \citep[version 2.4.0]{Paszke2019} with opt\_einsum \citep{Smith2018} for the unoptimized ordering of dimensions with MKL 2022.1 as backend.
The function \texttt{torch.einsum} results in similar performance to the unoptimized variant with Eigen and MKL.
This underlines that we choose the optimal contraction order and that the suggested optimization of the data layout speeds up the computation further compared to common tensor libraries such as pytorch.}

\subsection{Complete TT-AMEn algorithm}\label{subsec:performance_of_tt_amen}
For the complete algorithm for solving a linear system, we focus on the AMEn method as it needs at least an order of magnitude fewer operations than the other methods.
\begin{figure}
\subcaptionbox{\revised{Full variant (AMEn+SVD)}}{\includegraphics{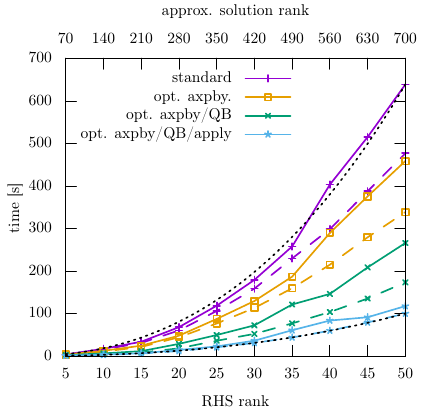}}
\hfill
\subcaptionbox{\revised{ALS variant (AMEn+ALS)}}{\includegraphics{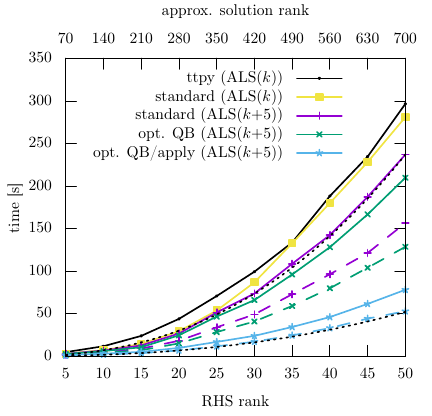}}
\caption{Timings for TT-AMEn for solving a linear system from a $50^{10}$ convection-diffusion problem ($c=10$) and random RHS $B_\text{TT}$ with varying ranks. Dashed lines use the TT-rank-1 preconditioner.
\revised{%
Dotted black lines illustrate the asymptotic complexity using the formula $c(0.35(r/700)^3+0.65(r/700)^2)$.
The heuristic ALS variant (right) is about twice as fast as the full variant (left).
For both variants, the time-to-solution is reduced by a factor of ${\sim}5$ by combining all suggested optimizations.
}}
\label{fig:tt_amen_performance}
\end{figure}
\revised{As shown in \autoref{fig:tt_amen_performance}, the full AMEn variant needs approximately twice the time of the ALS variant.
This is due to calculating the global residual $\mathcal A_\text{TT} X_\text{TT} - B_\text{TT}$ which is almost as costly as calculating $X_\text{TT}$ itself.
We also show timings obtained with the ttpy implementation (ALS variant) from \cite{Dolgov2014a} for which we linked with Intel MKL for the underlying operations.
Through using optimized building blocks and the suggested TT-rank-1 preconditioner, we speed up the calculation by a factor of ${\sim}5$ for both the full and the ALS variants.
In our tests with the ALS variant, we obtain better convergence and time-to-solution by using a slightly better approximation of the residual (cases with $\operatorname{ALS}(k+5)$ where $k$ denotes the AMEn enrichment rank and $k=\operatorname{rank}(B_\text{TT})$).
The TT-rank-1 preconditioner reduces the required number of inner GMRES iterations by about a factor of two here: from ${\sim}1530$ to ${\sim}750$ for the SVD variant and from ${\sim}1580$ to ${\sim}850$ for the ALS variant for the largest case.
However, the number of outer iterations (sweeps) stays the same and a significant part of the runtime is spent in the outer iteration.
Therefore, the preconditioner only speeds up the total runtime by about a factor of $1.2$--$1.5$.
}

\section{Conclusion and future work}\label{sec:conclusion}
In this paper, we discussed the complexity and the performance of linear solvers in tensor-train format.
In particular, we considered three different common methods, namely TT-GMRES, MALS (DMRG approach for linear systems) and AMEn, and tested their behavior for a simple, non-symmetric discretization of a convection-diffusion equation.
Concerning the complexity in terms of floating-point operations, we illustrated that AMEn can be about $100\times$ faster than MALS, which in turn can be about $100 \times$ faster than TT-GMRES.
These results already include an optimized orthogonalization scheme for the Arnoldi iteration in the TT-GMRES method which is also used as inner iteration of the MALS method.

Concerning the performance, we focussed on the required building blocks on a many-core CPU.
We suggested three improvements over the standard implementation: (a) exploiting orthogonalities in the TT-addition with subsequent truncation, (b) using a Q-less tall-skinny QR (TSQR) implementation to speed up costly singular value and QR decompositions, and (c) optimizing the memory layout/ordering of required tensor-contraction sequences for applying the tensor-train operator in the inner iteration.
As improvements (a) and (b) lead to less robust underlying linear algebra operations, we discussed their accuracy in the context of the required tensor-train operations.
In addition, we presented a simple generic preconditioner based on a tensor-train rank-1 approximation of the operator.
Overall, we obtained a speedup of about $5\times$ over the reference implementation on a 64-core CPU.

For future work, we want to investigate building block improvements for other tensor-network algorithms which are often based on very similar underlying operations.





\bibliographystyle{SageH}

\end{document}